\documentclass[11pt]{amsart}
\usepackage{amsfonts}
\usepackage{amssymb, amscd}
\usepackage{amsmath}
\input xy
\xyoption {all}

\usepackage{amsmath}
\usepackage{amsthm}
\usepackage{amssymb}
\usepackage{latexsym}

\newtheorem{lem}{Lemma}[section]
\newtheorem{cor}{Corollary}[section]

\newtheorem{thm}{Theorem}[section]
\newtheorem{defi}{Definition}[section]
\newcommand{\bb}[1]{\mbox{$\mathbb{#1}$}}
\theoremstyle{definition}

\theoremstyle{remark}
\newtheorem{rem}{Remark}[section]

\title{Specialization of motivic Hodge-Chern classes}

\begin{document}

\bibliographystyle{plain}

\author[J. Sch\"urmann ]{J\"org Sch\"urmann}
\address{J.  Sch\"urmann : Mathematische Institut,
          Universit\"at M\"unster,
          Einsteinstr. 62, 48149 M\"unster,
          Germany.}
\email {jschuerm@math.uni-muenster.de}

\begin{abstract}
In this paper we give a proof of the fact, that the motivic Hodge-Chern class  transformation
$MHC_y$ and Hirzebruch class transformation $MHT_{y*}$ for mixed Hodge modules and strictly specializable filtered 
${\mathcal D}$-modules commute with specialization in the algebraic and in a suitable complex analytic context.
Here specialization in the Hodge- and ${\mathcal D}$-module context means the corresponding
nearby cycles defined in terms of the $V$-filtration of Malgrange-Kashiwara. 
This generalizes a corresponding specialization result of Verdier about
MacPherson's Chern class transformation $c_*$.
\end{abstract}

\maketitle

\section{Introduction}
Let $k\subset \bb{C}$ be a subfield of the complex numbers. Using Saito's deep theory of {\em algebraic mixed Hodge modules} over $k$ (\cite{Sa1}-\cite{Sa7}),
we introduced in \cite{BSY} the {\em motivic Chern and Hirzebruch class transformations} 
$mC_y,MHC_y$ and $T_{y*},MHT_{y*}$ as natural transformations (cummuting with push down for proper morphisms)
fitting into a commutative diagram (with $\bb{Q}[y]_{loc}:=\bb{Q}[y,y^{-1},(1+y)^{-1}]$):
\begin{equation*}
\begin{CD}
 G_0(X)[y] @>>> G_0(X)[y,y^{-1}] @= G_0(X)[y,y^{-1}] \\
@A mC_y AA @A mC_y AA @AA MHC_y A \\
K_0(var/X) @>>> {\mathcal M}(var/X) @> \chi_{Hdg}>> K_0(MHM(X/k))\\
@V T_{y*} VV @V T_{y*} VV @VV MHT_{y*} V \\
H_*(X)\otimes \bb{Q}[y] @>>>  H_*(X)\otimes \bb{Q}[y,y^{-1}] @>>>  
H_*(X)\otimes \bb{Q}[y]_{loc} \:.
\end{CD}
\end{equation*}
Here $G_0(X)$ resp. $K_0(MHM(X/k))$ is the Grothendieck group of coherent sheaves resp. algebraic mixed Hodge modules on $X$,
$K_0(var/X)$ is the relative  Grothen\-dieck group af algebraic $k$-varieties over $X$, with ${\mathcal M}(var/X):=K_0(var/X)[\bb{L}^{-1}]$
its localization with respect to the class of the affine line 
$$\bb{L}=[\bb{A}^1_k\to pt]\in K_0(var/pt)$$ 
(compare e.g. \cite{Bi}). and finally
$H_*(X)=CH_{*}(X)$ is the Chow homology group \cite{Ful}, or for $k=\bb{C}$,  $H_*(X)=H_{2*}^{BM}(X)$ can also be the Borel-Moore homology
of $X$ (in even degrees).\\

The motivic Chern class transformations $mC_y, MHC_y$ are a $K$-theoretical refinement of the Hirzebruch class transformations $T_{y*},MHT_{y*}$,
since they are related by a (functorial) commutative diagram :
 \begin{equation*}
\begin{CD}
 {\mathcal M}(var/X) @> \chi_{Hdg}>> K_0(MHM(X/k)) @> MHC_y >>  G_0(X)[y,y^{-1}]\\
@V T_{y*} VV @V MHT_{y*} VV @VV td* V \\
H_*(X)\otimes \bb{Q}[y,y^{-1}] @>>>  H_*(X)\otimes\bb{Q}_{loc}  @< (1+y)^{-*}\cdot<<  H_*(X)\otimes \bb{Q}[y,y^{-1}] \:,
\end{CD}
\end{equation*}
with $td_*: G_0(X)\to H_*(X)\otimes \bb{Q}$ the {\em Todd class transformation} of Baum-Fulton-MacPherson \cite{BFM, Ful} and $(1+y)^{-*}\cdot$ the renormalization
given in degree $i$ by the multiplication
$$(1+y)^{-i}\cdot : H_i(-)\otimes\bb{Q}[y,y^{-1}] \to H_i(-)\otimes\bb{Q}[y,y^{-1},(1+y)^{-1}]
\:.$$
The characteristic class transformations $mC_y, T_{y*}$ are motivic refinements of the (rationalization of the) {\em Chern class transformation}
$$c_*: F(X)\to H_*(X)$$
of MacPherson \cite{M, Ken}, with $F(X)$ the abelian group of algebraically constructible functions, fitting into a (functorial) commutative diagram
\begin{equation*}
\begin{CD}
F(X) @< can << {\mathcal M}(var/X) \\
@V c_*\otimes \bb{Q} VV @VV T_{y*} V \\
H_*(X)\otimes \bb{Q} @ < y=-1<<  H_*(X)\otimes \bb{Q}[y,y^{-1}] \:.
\end{CD}
\end{equation*}

For $k=\bb{C}$, this was further improved in our recent survey \cite{Sch3}[Prop.5.2.1]. $MHT_{y*}$ factorizes as
$$MHT_{y*}: K_0(MHM(X)) \to H_*(X)\otimes \bb{Q}[y,y^{-1}]\subset  H_*(X)\otimes\bb{Q}[y]_{loc}\:,$$
fitting into a (functorial) commutative diagram
\begin{equation}
\begin{CD}
F(X) @<\chi_{stalk}<< K_0(D^b_c(X)) @< rat << K_0(MHM(X))\\
@V c_*\otimes \bb{Q} VV  @V c_*\otimes \bb{Q} VV @VV MHT_{y*} V \\
H_*(X)\otimes \bb{Q}  @= H_*(X)\otimes \bb{Q}  @ < y=-1<<  H_*(X)\otimes \bb{Q}[y,y^{-1}] \:.
\end{CD}
\end{equation}
Here $D^b_c(X)$ is the derived category of algebraically constructible sheaves (with rational coefficients) on $X$ (viewed as a complex analytic space), with
$$rat: D^bMHM(X) \to D^b_c(X) \quad \text{and} \quad rat: MHM(X) \to Perv(X)$$
associating to a (complex of) mixed Hodge module(s) the underlying (perverse) constructible sheaf complex.
Then 
$$K_0(D^b_c(X))=K_0(Perv(X))\quad \text{and}\quad K_0(D^bMHM(X))=K_0(MHM(X))$$ 
(see e.g. \cite{Sch2}[lem.3.3.1]), and $\chi_{stalk}$ is given by the Euler characteristic of the stalks.
Note that 
$$can: K_0(var/X)\to F(X) \quad \text{ and} \quad \chi_{stalk}\circ rat: K_0(MHM(X))\to F(X)$$ 
are surjective. Moreover, important functoriality results known for the MacPherson Chern class transformation $c_*$, like 
\begin{enumerate}
\item functoriality for proper morphisms,
\item multiplicativity for exterior products, and
\item a Verdier Riemann-Roch formula for smooth pullbacks,
\end{enumerate}
have been extended to $mC_y, T_{y*}$ in \cite{BSY} and to $MHC_y, MHT_{y*}$ in \cite{Sch3}.\\

The aim of this paper is to prove for these motivic transformations the counterpart of a famous result of Verdier,
that in the complex context the MacPherson Chern class transformation $c_*$ {\em commutes with specialization} \cite{V2}.
Let $f: X\to \bb{C}$ an algebraic function with $X_0:=\{f=0\}$.
Then Deligne's {\em nearby and vanishing cycle functors} (compare e.g. \cite{Sch2})
$$\Psi_f, \Phi_f: D^b_c(X) \to  D^b_c(X_0)$$
induce similar transformations for constructible functions fitting into a commutative diagram
\begin{equation}\label{nearby}\begin{CD}
K_0(D^b_c(X)) @> \Psi_f, \Phi_f >>  K_0(D^b_c(X_0)) \\
@V \chi_{stalk} VV @V \chi_{stalk} VV \\
F(X) @> \Psi_f, \Phi_f >> F(X_0)  \:.
\end{CD}\end{equation}
See also \cite{Sch} for an elementary description (without using sheaf theory) of $\Psi_f, \Phi_f: F(X)\to F(X_0)$ in terms of local {\em Milnor fibrations}.
Assume now that $X_0$ is a hypersurface of codimension one, so that one also has a homological {\em Gysin homomorphism} for the inclusion 
$i: X_0=\{f=0\}\to X$ (\cite{V2, Ful}):
$$i^!: H_*(X) \to H_{*-1}(X_0) \quad \text{and} \quad i^!: G_0(X)\to G_0(X_0)\:.$$
Then Verdier's specialization result can be formulated as the equality of the following two transformations:
 \begin{equation}\label{special-verdier}\begin{CD}
F(X) @> c_*\circ \Psi_f = > i^! c_* > H_*(X_0) \:.
\end{CD}\end{equation}
He also proved this result in the {\em complex analytic context for compact spaces}. Another proof in the {\em algebraic context over a base field $k$ of characteristic zero} was  
later given in \cite{Ken2}.\\

Consider now a base field $k\subset \bb{C}$, with $f: X\to \bb{A}^1_k$ a morphism and $X_0:=\{f=0\}$ the fiber over $0$.
Then one can consider the nearby- and vanishing cycle functors $\Psi_f$ and $\Phi_f$ either on the motivic level of localized relative Grothendieck groups 
$${\mathcal M}(var/-)=K_0(var/-)[\bb{L}^{-1}]$$ 
(see \cite{Bit,GLM}), or on the Hodge theoretical level of algebraic mixed Hodge modules over $k$
(\cite{Sa1, Sa2, Sa5}), ``lifting'' under the base change $k\to \bb{C}$ the corresponding functors on the level of  algebraically constructible sheaves 
(\cite{Sch2}) and  functions (\cite{Sch,V2}), so that the following diagram commutes:
\begin{equation}\label{all-nearby}\begin{CD}
K_0(var/X)\\
@VVV\\
{\mathcal M}(var/X) @> \Psi_f^m, \Phi_f^m >>  {\mathcal M}(var/X_0)\\
@V \chi_{Hdg} VV  @V \chi_{Hdg} VV  \\
K_0(MHM(X/k)) @> \Psi'^H_f, \Phi'^H_f >>  K_0(MHM(X_0/k))\\
@V rat VV @V rat VV \\
K_0\left(D^b_c(X(\bb{C}))\right) @> \Psi_f, \Phi_f >>  K_0\left(D^b_c(X_0(\bb{C}))\right) \\
@V \chi_{stalk} VV @V \chi_{stalk} VV \\
F(X(\bb{C})) @> \Psi_f, \Phi_f >> F(X_0(\bb{C}))  \:.
\end{CD}\end{equation}
Here we use the notation $\Psi'^H_f:= \Psi^H_f[1]$ and $\Phi'^H_f:=\Phi^H_f[1]$ for the shifted
functors, with 
$$\Psi^H_f, \Phi^H_f: MHM(X/k)\to MHM(X_0/k)$$ and 
$$\Psi_f[-1],\Phi_f[-1]: Perv(X(\bb{C}))\to Perv(X_0(\bb{C}))$$ 
preserving mixed Hodge modules and perverse sheaves, respectively. 
On the level of Grothendieck groups of mixed Hodge modules one has the relation
\begin{equation}
\Psi'^H_f = i^* + \Phi'^H_f: K_0(MHM(X/k)\to  K_0(MHM(X_0/k))\:,
\end{equation}
and similarly for the other transformations. So the vanishing cycles measure the 
difference between the pullback $i^*$ and the nearby cycles. In particular
$[\Psi'^H_f({\mathcal M}]= [i^*{\mathcal M}]$ in case $\Phi'^H_f({\mathcal M})=0$.

\begin{rem}  The motivic nearby and vanishing cycles functors of \cite{Bit,GLM}
take values in a refined {\em equivariant} localized Grothendieck group
${\mathcal M}^{\hat{\mu}}(var/X)$ of equivariant algebraic varieties over $X$ with a ``good'' action
of the pro-finite group $\hat{\mu}=\lim \mu_n$ of {\em roots of unity}. 
In our applications above we don't need to take this action into account.

Also note that for the commutativity of diagram (\ref{all-nearby}) one has to use $\Psi'^H_f$ (as opposed to 
$\Psi^H_f$, as stated in \cite{GLM}[Prop.3.17]; this fits in fact with the reference given in
the proof of loc.cit). Moreover, the Grothendieck group  ${\mathcal M}^{\hat{\mu}}(var/X)$ used in \cite{GLM}
is finer than the one used in \cite{Bit}. But both definitions of the motivic nearby and vanishing
cycle functors are compatible (\cite{GLM}[Rem.3.13]), and $\chi_{Hdg}$ also factorizes over 
${\mathcal M}^{\hat{\mu}}(var/X)$ in sense of \cite{Bit} by the same argument as for \cite{GLM}[(3.16.2)].
\end{rem}

Before formulating our main result, let us explain a motivating example in the case
$k=\bb{C}$, with $i: M_0=\{f=0\}\to M$ a codimension one inclusion of complex algebraic manifolds (with $m:=dim(M)$) and 
$${\mathcal M} \in MHM(M)[m]\subset D^bMHM(M)$$ 
corresponding to a {\em ``good'' variation 
of mixed Hodge structures} on $M$, i.e. such that $rat({\mathcal M})$ is a local system.
Then by \cite{Sch3} there is a cohomological characteristic class 
$MHC^y({\mathcal M})\in K^0(M)[y,y^{-1}]$ commuting with the pullback $i^*$, with $K^0(M)$ the Grothendieck group of algebraic
vector bundles or coherent locally free sheaves (of sections), such that
$$MHC_y([{\mathcal M}])= MHC^y({\mathcal M})\cap \left(\Lambda_y([\Omega^1_M])
\cap [{\mathcal O_M}]\right) \:.$$
Here the pairing
$$\cap=\otimes:  K^0(M)[y,y^{-1}]\times  G_0(M)[y,y^{-1}]\to  G_0(M)[y,y^{-1}]$$
is induced from the tensor product, with
$$\Lambda_y([T^*M]):=\Lambda_y([\Omega^1_M]):= \sum_{i\geq 0}\; [\Omega^i_M]\cdot y^i
\in K^0(M)[y] $$
the {\em total Lambda class} of the cotangent bundle $T^*M$. Then one gets
\begin{equation}\begin{split}
i^!MHC_y([{\mathcal M}])&=i^*\left(MHC^y({\mathcal M})\cup \Lambda_y([\Omega^1_M])\right)
\cap i^!([{\mathcal O_M}])\\
&=\left( MHC^y(i^*{\mathcal M})\cup \Lambda_y([i^*T^*M])\right) \cap [{\mathcal O_{M_0}}]\\
&= \Lambda_y([N^*_{M_0}M]) \cap MHC_y([i^*{\mathcal M}])\\
&= (1+y)\cdot MHC_y(i^*[{\mathcal M}])\:.
\end{split}\end{equation}
Here one uses the {\em multiplicativity} of $\Lambda_y(-)$ with respect to the short exact sequence of vector bundles
$$0\to N^*_{M_0}M \to T^*M|M_0 \to T^*M_0\to 0 \:,$$
and the triviality of the conormal bundle $N^*_{M_0}M$ (coming from the section $df$) so that
$$\Lambda_y(N^*_{M_0}M )= (1+y)\cdot [{\mathcal O_{M_0}}]\in K^0(M_0)[y]\:.$$
But in this special case there are no vanishing cycles $\Phi'^H_f({\mathcal M})=0$ so that
$$MHC_y(i^*[{\mathcal M}])= MHC_y\left(\Psi'^H_f([{\mathcal M}])\right)$$
and
$$i^!MHC_y([{\mathcal M}])=(1+y)\cdot MHC_y\left(\Psi'^H_f([{\mathcal M}])\right)\:.$$

And this formula holds in general, i.e. assume now that $X_0=\{f=0\}$ is a hypersurface of codimension one in the algebraic variety $X$ over $k\subset \bb{C}$. Then we prove the following counterpart of Verdier's
specialization result: 

\begin{thm} \label{motthm}
The motivic Hodge-Chern class transformation $MHC_y$ commutes with specialization in the following sense:
\begin{equation}
(1+y)\cdot MHC_y(\; \Psi'^H_f (-)\;) = 
-(1+y)\cdot MHC_y(\; \Psi^H_f (-)\;) = i^!MHC_y(-)   
\end{equation} 
as transformations $K_0(MHM(X/k))\to G_0(X_0)[y,y^{-1}]$.
\end{thm} 

\begin{rem}
 Our proof also gives the same result in the {\em embedded complex analytic} context,  with
$X\subset M$ a complex analytic subspace of the complex manifold $M$ and $f: M\to \bb{C}$ a holomorphic function
such that $X_0:=X\cap \{f=0\}$ is a {\em compact} hypersurface of codimesion one in $X$.
Here we use analytic mixed Hodge modules $MHM(X)$ on $X$ (\cite{Sa1, Sa2}), which we identify with analytic mixed Hodge modules 
$MHM_X(M)$ on $M$ with support in $X$ (and similarly for $X_0$). Here $G_0(X)$ is of course the Grothendieck group
of analytic coherent ${\mathcal O}_X$-sheaves, which has a canonical map to the Grothendieck group $G_0^X(M)$ of analytic coherent ${\mathcal O}_M$-sheaves
on $M$ with support in $X$. And $G_0(X_0)=G_0^{X_0}(M)$, since $X_0$ is compact.
\end{rem}

Let us come back to the algebraic context. Then another earlier result of Verdier \cite{V1, Ful} states that the Todd class transformation
of Baum-Fulton-MacPherson \cite{BFM, Ful}
$$td_*: G_0(-) \to H_*(-)\otimes \bb{Q}$$
{\em commutes} with the Gysin homomorphisms $i^!$ 
in these homology theories. The {\em unnormalized} Hirzebruch class transformation $MH\tilde{T}_{y*}$ of \cite{BSY} is defined as
$$MH\tilde{T}_{y*}:=td_* \circ MHC_y: K_0(MHM(-/k)) \to H_*(-)\otimes\bb{Q}[y,y^{-1}] \:.$$
Since $td_*$ commutes with specialization, we get the

\begin{cor} Consider an algebraic function $f: X\to \bb{A}^1_k$ such that the inclusion
$i: X_0:=\{f=0\}\to X$ of the zero fiber is everywhere of codimension one.
Then the Hirzebruch class transformations  $MH\tilde{T}_{y*}$ and $MHT_{y*}$ commute with
specialization in the following sense:
\begin{equation}
(1+y)\cdot MH\tilde{T}_{y*}(\; \Psi'^H_f (-)\;) = i^!MH\tilde{T}_{y*}(-)   
\end{equation} 
as transformations $K_0(MHM(X/k))\to H_*(X_0)\otimes\bb{Q}[y,y^{-1}]$.
\begin{equation}
MHT_{y*}(\; \Psi'^H_f (-)\;) = i^!MHT_{y*}(-)   
\end{equation} 
as transformations $K_0(MHM(X/k))\to H_*(X_0)\otimes\bb{Q}[y,y^{-1},(1+y)^{-1}]$. \hfill $\Box$
\end{cor}

By the definition of $\Psi_f^m$ in \cite{Bit,GLM} one has
$$\Psi_f^m(K_0(var/M))\subset im(K_0(var/X)\to {\mathcal M}(X))\:,$$
so that 
$${MHT_{y}}_*\circ \Psi_f^m : K_0(var/M)  \to
H_*(X)\otimes\bb{Q}[y]\subset H_*(X)\otimes\bb{Q}[y,y^{-1}] \:.$$
Together with \cite{Sch3}[Prop.5.2.1] one therefore gets the following commutative diagram
of specialization results, where one uses the base change $k\to \bb{C}$ for the last part:
\begin{equation}\label{special}\begin{CD}
K_0(var(X)) @> {T_{y}}_*\circ \Psi_f^m = > i^!\circ {T_{y}}_* > H_*(X_0)\otimes \bb{Q}[y]\\
@V \chi_{Hdg} VV @VVV\\
K_0(MHM(X/k)) @> {MHT_{y}}_*\circ \Psi'^H_f = > i^!\circ {MHT_{y}}_* > H_*(X_0)\otimes \bb{Q}[y,y^{-1}]\\
@V \chi_{stalk}\circ rat VV @VV y=-1 V \\
F(X(\bb{C})) @> c_*\circ \Psi_f = > i^! c_* > H_*(X_0(\bb{C}))\otimes \bb{Q} \:.
\end{CD}\end{equation}

Applications of this corollary to the calculation of the Hirzebruch class $$T_{y*}(X):=T_{y*}([id_X])$$
of a global hypersurface $X:=\{f=0\}$ inside a manifold $M$ in the complex algebraic context are
given in our recent joint paper \cite{CMSS}.\\

Note that the motivic Chern class transformation $MHC_y$ is defined 
in terms of the {\em filtered de Rham complex} of the filtered ${\mathcal D}$-module underlying a mixed Hodge module
(\cite{BSY, Sch3}). So in the embedded context $X\subset M$ a closed subspace of the manifold $M$,
with $f: M\to \bb{A}^1_k$, our main result theorem \ref{motthm} becomes a purely
${\mathcal D}$-module theoretic result about coherent ${\mathcal D_M}$-modules with a good (Hodge) filtration $F$ supported on $X$,
which are {\em strictly specializable} with respect to the hypersurface $\{f=0\}$.
Here one uses the ${\mathcal D}$-module description of nearby cycles in terms of the {\em $V$-filtration of Malgrange-Kashiwara}
(\cite{Mal, Ka, MM}). Then the strict specializability is a condition on the relation between the Hodge filtration $F$ and the $V$-filtration
introduced by M. Saito \cite{Sa1} in his definition of {\em pure Hodge modules}. So the filtered ${\mathcal D}$-module underlying a pure (and also a
mixed) Hodge module satisfies this technical condition (more or less) by definition.
And in the next section we will show that in the embedded context exactly this {\em  strict specializability} is needed for the proof
of our main result! Using {\em resolution of singularities and functoriality} under proper morphisms, we deduce in the last section our main result for
mixed Hodge modules on a singular space $X$ from the corresponding   ${\mathcal D}$-module result in the embedded context.

\section{Specialization for ${\mathcal D}$-modules}
We are working either in the algebraic context over a base field $k\subset \bb{C}$, i.e. all sheaves are algebraic sheaves in the Zariski topology,
or in the complex analytic context. Consider an algebraic (or analytic) manifold $M$ of dimension $n$, together with 
$$t=pr_2: M':=M\times \bb{A}^1_k \to \bb{A}^1_k$$
the projection onto the second factor. Of course in the complex analytic context, we consider
the affine line $\bb{C}=\bb{A}^1_C$ just as a complex manifold.
Let $i': M\simeq M\times\{0\} \to M':=M\times \bb{C}$
be the closed inclusion of the zero fiber. \\

Let ${\mathcal D}_{M'}$ be the corresponding coherent sheaf of algebraic or analytic differential operators. This has an increasing filtration $F$ by locally free coherent
${\mathcal O_{M'}}$-sheaves $F^i{\mathcal D}_{M'}$ ($i\in \bb{N}_0$) given by the differential operators of order less or equal $i$, with 
$$\cup_{i\geq 0}\: F^i{\mathcal D}_{M'}= {\mathcal D}_{M'}\:.$$
(see e.g. \cite{MM,Sa6}). 
Let ${\mathcal I}\subset {\mathcal O}_{M'}$ be the ideal sheaf defining $M$, i.e. the sheaf
of functions vanishing along $M$. Then
the {\em increasing $V$-filtration of Malgrange-Kashiwara} with respect to the smooth hypersurface $M\subset M'$ is  for $k\in \bb{Z}$ {\em locally} defined by
$$V_k{\mathcal D}_{M'}:=\{P\in {\mathcal D}_{M'}|\; P({\mathcal I}^{j+k})\subset {\mathcal I}^j
\quad \text{for all $j\in \bb{Z}$}\}\:.$$
Here ${\mathcal I}^j:={\mathcal O_{M'}}$ for $j<0$. Note that
$$\cap_{k\in \bb{Z}}\:V_k{\mathcal D}_{M'}=\{0\} \quad \text{and} \quad
\cup_{k\in \bb{Z}}\:V_k{\mathcal D}_{M'}={\mathcal D}_{M'}\:.$$
Moreover
$V_k{\mathcal D}_{M'}|\{f\neq 0\}= {\mathcal D}_{M'}|\{f\neq 0\}$ for all $k\in \bb{Z}$
so that $gr^V_k{\mathcal D}_{M'}$ is supported on $M$. By definition one has
$$t\in  V_{-1}{\mathcal D}_{M'}, \:\partial_t \in  V_{1}{\mathcal D}_{M'}
\quad \text{ and} \quad
\partial_t t = 1+t\partial_t \in  V_{0}{\mathcal D}_{M'} \:.$$ 
Similarly we have for the sheaf 
${\mathcal D}_{M'/A}$
of {\em relative differential operators along the fibers of $t$}: 
$${\mathcal D}_{M'/A}\subset  V_{0}{\mathcal D}_{M'}
\quad \text{and} \quad gr^V_0{\mathcal D}_{M'}|M={\mathcal D}_{M}[\partial_t t] \:.$$

\begin{defi}[Specializability]
A {\em coherent left}  ${\mathcal D}_{M'}$-module ${\mathcal M}$ is said to be
{\em specializable along $M$} if for any (closed) point $x\in M$ (with residue field $k_x$) and any local section
$m\in {\mathcal M}_x$ there is a nonzero polynomial $b(s)\in k_x[s]$ such that
\begin{equation}\label{b-fct}
b(\partial_t t)\cdot m \in V_{-1}{\mathcal D}_{M',x}\cdot m \:.
\end{equation}
The {\em Bernstein polynomial (or $b$-function)} $b_m$ of $m\in  {\mathcal M}_x$
is the monic polynomial in $k_x[s]$ of minimal degree satisfying the condition (\ref{b-fct}).
\end{defi}

For example any {\em holonomic} ${\mathcal D}_{M'}$-module is specializable along $M$
(compare e.g. \cite{MM}[prop.4.4.2]). Similarly, for a short exact sequence
$$0\to  {\mathcal M}'\to  {\mathcal M}\to  {\mathcal M}''\to 0$$
of coherent  ${\mathcal D}_{M'}$-modules one has ${\mathcal M}$ is specializable
if and only if  ${\mathcal M}'$ and  ${\mathcal M}''$ are specializable (\cite{MM}[prop.4.2.4]).\\

Let us  fix a total order $\leq$ on $\bb{C}$ with the following properties for all
$\alpha,\beta \in \bb{C}$ and $a\in \bb{R}$:
\begin{enumerate}
\item[(o1)] it induces the usual order on $\bb{R}$,
\item[(o2)] $\alpha +a < \beta + a \Leftrightarrow \alpha < \beta$,
\item[(o3)]$\alpha < \alpha +1$ and there is some $m\in\bb{N}$ with $\alpha< \beta+m$.
\end{enumerate}

\begin{defi}[Canonical $V$-filtration]\label{can-V}
Let ${\mathcal M}$ be a coherent {\em left} ${\mathcal D}_{M'}$-module, which is specializable with respect to the hypersurface $M$. Then the {\em canonical increasing} $V$-filtration
of ${\mathcal M}$ (indexed by $\bb{C}$ and $<$) is defined as
\begin{equation}
V_{\alpha}{\mathcal M}_x:=\{m\in {\mathcal M}_x|\:\text{all roots of the $b$-function $b_m$ are $\geq -\alpha -1$}\:\}
\end{equation}
for all (closed) points $x\in M$ (compare \cite{MM}[def.4.3.3]).
In the complex analytic context we assume in addition that $M\cap supp({\mathcal M})$ is {\em compact}.
Then the canonical $V$-filtration is indexed {\em  discretely} by $A+\bb{Z} \subset 
\bar{k} \subset\bb{C}$ for a finite subset $A\subset \bar{k} \subset\bb{C}$, and it is the 
{\em unique} such filtration $V$
with the following properties:
\begin{enumerate}
\item $\bigcup_{\alpha} V_{\alpha} {\mathcal M} = {\mathcal M}$, and each $V_{\alpha} {\mathcal M}$ is a coherent
$V_{0}{\mathcal D}_{M'}$ module.
\item $(V_k{\mathcal D}_{M'}) (V_{\alpha} {\mathcal M}) \subset V_{\alpha+k} {\mathcal M}$ for all $\alpha\in \bb{C}, k\in \bb{Z}$, \\
and $t(V_{\alpha} {\mathcal M}) = V_{\alpha-1} {\mathcal M}$ for all $\alpha <0$.
\item $\partial_t t+\alpha$ is {\em nilpotent} on $gr^V_{\alpha}{\mathcal M}:= 
V_{\alpha} {\mathcal M}/V_{<\alpha} {\mathcal M}$,\\
 with $V_{<\alpha} {\mathcal M}:=V_{\beta} {\mathcal M}$
for $\beta:=\sup \{\beta'\in A+\bb{Z}|\; \beta'< \alpha\}$.
\end{enumerate}
\end{defi}

These conditions imply $t: gr^V_{\alpha}{\mathcal M}\to gr^V_{\alpha-1}{\mathcal M}$ is bijective for all
$\alpha \neq 0$, and $\partial_t: gr^V_{\alpha}{\mathcal M}\to gr^V_{\alpha+1}{\mathcal M}$ is bijective for all
$\alpha \neq -1$. 
Moreover (compare \cite{MM}[prop.4.4.3] and \cite{Sa1}[lem.3.1.4])
\begin{equation} \label{injective}
t\cdot: V_{\alpha} {\mathcal M}\to V_{\alpha-1} {\mathcal M}\quad
\text{is injective for all $\alpha<0$.}
\end{equation}
Finally all $gr^V_{\alpha}{\mathcal M}|M$ are coherent {\em left} ${\mathcal D}_{M}$-modules (\cite{MM}[cor.4.3.10]), which are also all (regular) holonomic, if ${\mathcal M}$ has this property (\cite{MM}[cor.4.6.3, cor.4.7.5]). In the algebraic context
over a base field $k$ this means by definition that the corresponding 
${\mathcal D}$-modules over $M'(\bb{C})$ and $M(\bb{C})$ obtained by the base extension
$k\to \bb{C}$ are (regular) holonomic.

\begin{defi}
Let ${\mathcal M}$ be a coherent {\em left} ${\mathcal D}_{M'}$-module, which is specializable with respect to the hypersurface $M$. Then ${\mathcal M}$ is called {\em quasi-unipotent}
(along $M$), if for all $m\in {\mathcal M}_x$ with $x$ a (closed) point of $M$ the Berstein
polynomial $b_m(s)\in k_x[s]$ has only rational roots in $\bb{Q}\subset k\subset \bar{k}$.
Or in other words, if the canonical $V$-filtration is discretely indexed by $\bb{Q}$
(so that the total order $<$ on $\bb{C}$ above isn't needed).
\end{defi}

\begin{rem}
Here we are working with an {\em increasing} $V$-filtration for {\em left} ${\mathcal D}$-modules,
so that $\pm (\partial_t t+\alpha)$ is {\em nilpotent} on $gr^V_{\alpha}{\mathcal M}$.
\begin{enumerate}
\item If we work with the associated {\em analytic}  ${\mathcal D}$-module ${\mathcal M}^{an}$, then
$(V_{\alpha}{\mathcal M})^{an}$ is the corresponing $V$-filtration in the analytic context
(by uniqueness). Similarly in the algebraic context over a base field $k\subset \bb{C}$,
if one takes the base extension $k\to \bar{k}$ or $k\to \bb{C}$.
\item If we switch to the corresponding {\em right} ${\mathcal D}$-module 
${\mathcal M}^r :=\omega_{M'}\otimes {\mathcal M}$,
with the induced $V$-filtration $\omega_{M'}\otimes V_{\alpha}{\mathcal M}$, then
$t\partial_t -\alpha$ is nilpotent on $gr^V_{\alpha}{\mathcal M}^r$
(fitting with the convention of \cite{Sa1}[Def. 3.1.1]).
\item If we use the corresponding {\em decreasing} $V$-filtration $V^{\alpha}:=V_{-\alpha-1}$
for {\em left} ${\mathcal D}$-modules, then $t\partial_t -\alpha$ is nilpotent on $gr_V^{\alpha}{\mathcal M}$
(fitting with the convention of \cite{Sa1}[Introduction, p.851]).
\end{enumerate}
\end{rem}

Assume now that ${\mathcal M}$ in addition is endowed with a {\em good filtration $F$},
i.e. with an increasing filtration $F_p{\mathcal M}$ ($p\in \bb{Z}$) by coherent
${\mathcal O}_{M'}$-modules, such that 
\begin{enumerate}
\item $F_p{\mathcal M}=0$ for $p<<0$ small enough,
\item $\cup_{p\in \bb{Z}}\; F_p{\mathcal M}= {\mathcal M}$,
\item $(F_k{\mathcal D}_{M}) (F_p {\mathcal M}) \subset F_{k+p} {\mathcal M}$ for all $k,p\in \bb{Z}$,
\item $gr^F{\mathcal M}$ is a coherent $gr^F{\mathcal D}_{M'}$ module. 
\end{enumerate}
Then one can ask about the relation between the two filtrations $F$ and $V$.
The following notion will be important for our main results (compare \cite{Sa1}[Def.3.2.1, p.905] and \cite{Sab2}[conditions 4.3.(b-c), p.46]): 

\begin{defi}[Strict specializability]
Let ${\mathcal M}$ be a coherent {\em left} ${\mathcal D}_{M'}$-module, which is specializable with respect to the hypersurface $M$. Assume that it is also endowed with a good filtration $F$.
The filtered ${\mathcal D}_{M'}$-module $({\mathcal M},F)$ is called {\em strictly specializable} along $M$, iff
\begin{enumerate}
\item[s1.]  The induced $F$-filtration on $gr^V_{\alpha} {\mathcal M}$ is good for all $-1\leq \alpha \leq 0$.
\item[s2.] $t(F_p V_{\alpha} {\mathcal M}) = F_p V_{\alpha-1} {\mathcal M}$ for all $\alpha <0$
and $p\in \bb{Z}$.
\item[s3.] $\partial_t (F_p gr^V_{\alpha} {\mathcal M}) = F_{p+1} gr^V_{\alpha+1} {\mathcal M}$ for all $\alpha >-1$ and $p\in \bb{Z}$.
\end{enumerate}
\end{defi}
Here the {\em induced} $F$-filtration on $V_{\alpha}{\mathcal M}$ resp. $gr^V_{\alpha} {\mathcal M}$ are given by (cf. \cite{De})
$$F_pV_{\alpha}{\mathcal M}:= F_p{\mathcal M}\cap  V_{\alpha}{\mathcal M}$$ 
resp.
$$F_p gr^V_{\alpha} {\mathcal M}:=im\left(F_pV_{\alpha}{\mathcal M}\to gr^V_{\alpha} {\mathcal M}\right) \simeq F_pV_{\alpha}{\mathcal M}/ F_pV_{<\alpha}{\mathcal M} \:.$$
Then the induced short exact sequence of filtered ${\mathcal O}_{M'}$-modules
\begin{equation}
0\to ( V_{<\alpha}{\mathcal M},F) \to ( V_{\alpha}{\mathcal M},F) \to 
(gr^V_{\alpha} {\mathcal M},F)\to 0
\end{equation}
is {\em strict} (compare \cite{De, La}) in the sense that
$$0\to F_pV_{<\alpha}{\mathcal M} \to F_pV_{\alpha}{\mathcal M} \to 
F_pgr^V_{\alpha} {\mathcal M}\to 0$$
is exact for all $p\in \bb{Z}$, which also implies that
$$0\to gr^F_pV_{<\alpha}{\mathcal M} \to gr^F_pV_{\alpha}{\mathcal M} \to 
gr^F_pgr^V_{\alpha} {\mathcal M}\to 0$$
is exact for all $p\in \bb{Z}$. Finally one has a {\em canonical isomorphism}
(compare \cite{De}[(1.2.1)]):
\begin{equation}
gr^F_pgr^V_{\alpha} {\mathcal M}\simeq gr^V_{\alpha}gr^F_p {\mathcal M}
\quad \text{for all $p\in \bb{Z}, \alpha \in \bb{C}$.}
\end{equation}
And similarly for $V_{\alpha}/V_{\beta}$ instead of $gr^V_{\alpha}$ for $\alpha<\beta$.

\begin{rem} \label{pure}
By Saito's work \cite{Sa1, Sa2, Sa5}, the underlying filtered $D$-module of an {\em algebraic or analytic pure Hodge module} on $M'$ is {\em strictly specializable and  quasi-unipotent} along $M$.
\begin{enumerate}
\item In the analytic context this is so by definition (where we again assume that the intersection of the support with $M$ is compact, so that the canonical $V$-filtration is indexed discretely).
\item In the complex algebraic context a  pure Hodge module of \cite{Sa2} is by definition ``extendable and 
quasi-unipotent'' at infinity so that one can assume $M$ is compact and then the properties
(s1-s3) in the analytic context of \cite{Sa1} imply by GAGA and flatness of
${\mathcal O}^{an}$ over ${\mathcal O}$ the same properties for the underlying algebraic filtrations.
\item In the algebraic context over a base field $k\subset \bb{C}$ the claim follows from the
complex algebraic context by the exact and flat base change $k\to \bb{C}$ (compare \cite{Sa5}[rem. on p.9]).
\item Then also the underlying filtered $D$-module of an {\em algebraic or analytic mixed Hodge module} on $M'$ as in \cite{Sa2, Sa5} is {\em strictly specializable and  quasi-unipotent} along $M$, since it is a finite successive extension of pure Hodge modules (by the weight filtration). And the corresponding short exact extension sequences of mixed Hodge modules
a strict with respect to the (Hodge) filtration $F$. 
Since also the canonical $V$-filtration behaves well under extensions, one easily gets the
quasi-unipotence and the condition (s1-s3) above by induction (using the properties of the underlying $V$-filtrations from definition \ref{can-V}).
\end{enumerate}
\end{rem} 

For a closed algebraic (or analytic) subset $X'\subset M'$ we get a cartesian diagram of inclusions:
$$\begin{CD}
X:=X\cap \{t=0\} @>i >> X'\\
@VVV @VVV\\
M @>> i' > M'\:. 
\end{CD}$$
 
And we are interested in the algebraic (or analytic) {\em de Rham complex} $DR^*({\mathcal M})$ of a filtered ${\mathcal D}_{M'}$-module $({\mathcal M},F)$ supported on $X'$:
\begin{equation}
\begin{CD}
[ \cdots \Omega^k_{M'} \otimes {\mathcal M} @> \nabla^k >>   \Omega^{k+1}_{M'} \otimes {\mathcal M} \cdots]\:, 
  \end{CD}
\end{equation}
with $\Omega^{n+1}_{M'} \otimes {\mathcal M}$ put in {\em degree zero}
(and $n=dim(M)$, $n+1=dim(M')$).
Here $\nabla^k$ comes from the {\em integrable connection} $\nabla$ on ${\mathcal M}$ given by the {\em left} ${\mathcal D}_{M'}$-module structure:
\begin{equation}
\nabla^k(\omega \otimes m) = d\omega \otimes m + (-1)^k \omega \wedge \nabla m \:.
\end{equation}
From this it follows that the de Rham complex gets an induced filtration 
\begin{equation}
\begin{CD} F_p DR^*({\mathcal M}) =
[ \cdots \Omega^k_{M'} \otimes F_{p+k}{\mathcal M} @> \nabla^k >>   \Omega^{k+1}_{M'} \otimes F_{p+k+1}{\mathcal M} \cdots]\:, 
  \end{CD}
\end{equation}  
such that differentials $gr^F_p(\nabla^k)$ of the graded complex
\begin{equation}
\begin{CD} gr^F_p DR^*({\mathcal M}) =
[ \cdots \Omega^k_{M'} \otimes gr^F_{p+k}{\mathcal M} @> gr^F_p(\nabla^k) >>   \Omega^{k+1}_{M'} \otimes gr^F_{p+k+1}{\mathcal M} \cdots] 
  \end{CD}
 \end{equation} 
are ${\mathcal O}_{M'}$-linear! Since the filtration $F$ of ${\mathcal M}$ is {\em good}, one also knows (compare \cite{Sa1}[Lem.2.1.17, p.882]) that
\begin{equation}
gr^F_p DR^*({\mathcal M}) \in D^b_{coh,X'}(M')
\end{equation}
is a bounded complex with coherent cohomology supported on $X'$, which is acyclic for almost all $p$.

\begin{defi}\label{MHC}
The {\em motivic Hodge-Chern class} $MHC_y(({\mathcal M},F))$ of the coherent ${\mathcal D}$-module 
$({\mathcal M},F)$ with its good filtration $F$ is given by
\begin{equation}
MHC_y(\;({\mathcal M},F)\;) := \sum_{p} [gr^F_{-p} DR^*({\mathcal M})]\cdot (-y)^p \in G_0(X')[y,y^{-1}] \:.
\end{equation}
\end{defi}
Here we first identify the derived Grothendieck group
$$K_0\left(D^b_{coh,X'}(M')\right)=G^{X'}_0(M')$$
with the Grothendieck group $G^{X'}_0(M')$ of coherent 
${\mathcal O}_{M'}$-sheaves with support in $X'$, and the class $[-]$ of a complex means the alternating sum of the classes of its cohomology sheaves.
Next the closed inlusion $\kappa: X'\to M'$ induces a natural group homomorphism
on the Grothendieck group $G_0(X')$ of coherent 
${\mathcal O}_{X'}$-sheaves:
$$\kappa_*: G_0(X')\to G^{X'}_0(M')\:.$$
By the Hilbert-R\"{u}ckert theorem this is an isomorphism in the algebraic context,
as well as in the analytic context for $X'$ compact (or relatively compact,
if we are ``allowed to take a shrinking'').

\begin{rem}
If  $({\mathcal M},F)$ underlies a mixed Hodge module with support in $X'$,
then all $gr^F_p{\mathcal M}$ are in fact ${\mathcal O}_{X'}$-modules so that
$gr^F_p DR^*({\mathcal M})$ is already a well defined complex in $D^b_{coh}(X')$!
Just take locally a function on $M'$ vanishing along $X'$.
Then  $({\mathcal M},F)$ is also ``strictly specializable along $g$'', which implies
$$g(F_p{\mathcal M})\subset F_{p-1}{\mathcal M} \quad \text{and} \quad
g(gr^F_p{\mathcal M})=0$$
by \cite{Sa1}[lem.3.2.6].
\end{rem}

Now we want to compare the specialization (or nearby cycles)
$\Psi_t(({\mathcal M},F))$
of the {\em strictly specializable} pair $({\mathcal M},F)$ with the specialization
$i^!$ of the motivic Hodge-Chern classes, where the Gysin homomorphism 
$$i^!: G_0(X')\to G_0(X)$$
is induced by the derived pullback  map 
$$Li^*: D^b_{coh,X'}(M')\to D^b_{coh,X}(M)\:.$$
By the short exact sequence
$$\begin{CD} 0@>>> {\mathcal O}_{M}@>>>  {\mathcal O}_{M'} @> t\cdot >> {\mathcal O}_{M'}
@>>> 0 \:,
\end{CD}$$
this Gysin homomorphism $i^!$ is also represented by taking the (class of the)
tensor product with the complex $[t\cdot: {\mathcal O}_{M'}\to {\mathcal O}_{M'}]$.
In particular
\begin{equation}\label{i!i*}
i^!i_*=0: G_0(X)\to G_0(X); \: [{\mathcal F}]\mapsto [ {\mathcal F} \stackrel{t\cdot}{\to}
{\mathcal F}] \:.
\end{equation}

Let us now recall the 
\begin{defi} \label{psi}
The {\em nearby cycles} $\Psi_t(({\mathcal M},F))$ of the strictly specializable pair
$({\mathcal M},F)$ are given by
\begin{equation}
\Psi_t(({\mathcal M},F)):= \sum_{-1\leq \alpha <0} \Psi_{t,\alpha}(({\mathcal M},F)) \:,
\end{equation}
with $\Psi_{t,\alpha}(({\mathcal M},F)):=(gr^V_{\alpha}{\mathcal M}|M,F)$.\\
The {\em unipotent vanishing cycles} $\Phi_{t,uni}(({\mathcal M},F))$ of $({\mathcal M},F)$ are given by
\begin{equation}
\Phi_{t,uni}(({\mathcal M},F)):= (gr^V_{0}{\mathcal M}|M,F[-1]) \:,
\end{equation}
with the shifted filtration defined as $(F[k])_i:=F_{i-k}$.
\end{defi}

\begin{rem} Our defintion of the induced filtration fits with \cite{Sa1}[intro., p.851],
since we are using {\em left} ${\mathcal D}$-modules. For the corresponding {\em right}
${\mathcal D}$-modules one has to shift these induced $F$-filtrations by $[+1]$
(compare \cite{Sa1}[(5.1.3.3) on p.953]). This corrects then the different switching by
$$\otimes (\omega_{M'},F) \quad \text{ or}  \quad \otimes (\omega_{M},F)$$ 
from filtered left ${\mathcal D}$-modules
to filtered right ${\mathcal D}$-modules on $M'$ or $M$, with $F$ the trivial filtration
such that $gr^F_{-k}(-)= 0$ for $k\neq$ the dimension of the ambient manifold $M'$ or $M$.
\end{rem}
Now we can formulate our main result of this section.

\begin{thm} \label{mainthm}
Let $({\mathcal M},F)$ be a coherent ${\mathcal D}_{M'}$-module 
${\mathcal M}$ supported on $X'$ endowed with a good filtration $F$, 
which is strictly specializable along $M$. Then the motivic Hodge-Chern classes commute with specialization in the following sense:
\begin{equation}
-(1+y) \cdot MHC_y(\; \Psi_t(({\mathcal M},F))\;) = i^! MHC_y(\;({\mathcal M},F)\;) 
\end{equation}
as classes in $G_0(X)[y,y^{-1}]$.
\end{thm}

The proof is given on the next pages. It uses of course the properties (s1.)-(s3.). From  (s3.) one gets the
\begin{lem} $$\partial_t: ({\mathcal M}/V_{\alpha} {\mathcal M},F_*) \to ({\mathcal M}/V_{\alpha+1} {\mathcal M},F_{*+1})$$ 
is a {\em filtered isomorphism} for all $\alpha \geq -1$, so that also
$$\partial_t: gr^F_p({\mathcal M}/V_{\alpha} {\mathcal M}) \to gr^F_{p+1}({\mathcal M}/V_{\alpha+1} {\mathcal M})$$ 
is an {\em isomorphism} for all $\alpha \geq -1$ and $p\in \bb{Z}$.  \hfill $\Box$
\end{lem}
This follows by  $\bigcup_{\alpha} V_{\alpha} {\mathcal M} = {\mathcal M}$ from the corresponding result:
$$\partial_t: (V_{\beta}/V_{\alpha} {\mathcal M},F_*) \to (V_{\beta+1}/V_{\alpha+1} {\mathcal M},F_{*+1})$$ 
is a {\em filtered isomorphism} for all $\beta>\alpha \geq -1$.\\
And this follows by induction from (s3.) and the short exact sequence
$$\begin{CD}
0 \to (V_{\beta}/V_{\alpha} {\mathcal M},F_*) @>>> (V_{\beta'}/V_{\alpha} {\mathcal M},F_*) @>>> \cdots  \\
@V \:\:V \partial_t V  @VV \partial_t V \\
0 \to (V_{\beta+1}/V_{\alpha+1} {\mathcal M},F_{*+1}) @>>> (V_{\beta'+1}/V_{\alpha+1} {\mathcal M},F_{*+1}) @>>> \cdots
  \end{CD}$$
$$\begin{CD}
\cdots  (gr^V_{\beta'}{\mathcal M},F_*)  \to 0\\
@V\partial_t  V\:\: V \\
\cdots (gr^V_{\beta'+1}{\mathcal M},F_{*+1})  \to 0 \:,
\end{CD}$$
where $\beta' \in A+\bb{Z}$ is smallest number bigger than $\beta$.
Of course we also use the fact that $\partial_t: gr^V_{\beta}{\mathcal M}\to gr^V_{\beta+1}{\mathcal M}$ is bijective for all $\beta> -1$. \\

Similarly one gets from (s2.) and (\ref{injective})
 by induction over $p$ (with $F_p(-)=0$ for $p<<0$) the
\begin{lem} The multiplication
$$t\cdot : gr^F_p (V_{\alpha} {\mathcal M}) \to gr^F_p (V_{\alpha-1} {\mathcal M})$$
is an {\em isomorphism} for all $\alpha <0$ and $p\in \bb{Z}$. In particular
$$Li^* gr^F_p (V_{\alpha} {\mathcal M})= i^* gr^F_p (V_{\alpha} {\mathcal M})= 
gr^F_p (\;V_{\alpha}/ V_{\alpha-1}{\mathcal M}|M\;) $$
for all $\alpha <0$ and $p\in \bb{Z}$. \hfill $\Box$
\end{lem}

Now we can come to the main geometric idea. Since we work on the product manifold
$M'=M\times \bb{A}^1_k$ (which usually comes from a graph embedding for an algebraic (or analytic) function $f: M\to \bb{A}^1_k$), we have the splitting 
$$\Omega^1_{M'} = \Omega^1_{M'/A} \oplus  \Omega^1_{A}\:,$$
with $\Omega^1_{M'/A}$ the relative $1$-forms with respect to the submersion $t$
and $\Omega^1_{A} \simeq {\mathcal O}_{M'}dt$.
So we get for $0\leq k \leq n+1$:
$$\Omega^k_{M'} = \Omega^k_{M'/A} \oplus \Omega^1_{A}\otimes \Omega^{k-1}_{M'/A}\:.$$
And this induces a splitting of the de Rham complex $DR^*({\mathcal M})$ as a double complex
$$\begin{CD} DR^*_{/A}({\mathcal M}) @=
[ \cdots \Omega^k_{M'/A} \otimes {\mathcal M} @> \nabla_{/A}^k >>   \Omega^{k+1}_{M'/A} \otimes {\mathcal M} \cdots]\\
@V  \nabla_t VV @V \:\: \nabla_t^k VV @V  \nabla_t^{k+1} V \:\:V \\
\Omega^1_{A}\otimes DR^*_{/A}({\mathcal M}) @=
[ \cdots \Omega^1_{A}\otimes \Omega^k_{M'/A}\otimes {\mathcal M} @> \nabla_{/A}^k >>   \Omega^1_{A} \otimes\Omega^{k+1}_{M'/A}\otimes  {\mathcal M} \cdots] \:,
 \end{CD}$$
with the "top-dimensional forms" $\Omega^1_{A} \otimes\Omega^{n}_{M'/A}$ in bidegree
$(0,0)$.\\

Here the horizontal lines come from the corresponding {\em relative de Rham complex} 
$DR^*_{/A}({\mathcal M})$ of 
${\mathcal M}$ viewed only as left ${\mathcal D}_{M'/A}$-module, whereas the vertical maps are given by:
$$\nabla_t^k: \Omega^k_{M'/A} \otimes 
{\mathcal M}\to \Omega^1_{A} \otimes\Omega^k_{M'/A}\otimes  {\mathcal M}; $$
$$\omega \otimes m \mapsto dt\otimes (\partial_t\omega)\otimes m + 
(-1)^k dt\otimes\omega \otimes (\partial_t m) \:.$$

In particular $DR^*({\mathcal M})$  becomes a bifiltered double complex by
$$\begin{CD} F_pV_{\alpha}DR^*_{/A}({\mathcal M}) @= \\
@V  \nabla_t VV\\
\Omega^1_{C}\otimes F_{p+1}V_{\alpha+1}DR^*_{/A}({\mathcal M}) @=
  \end{CD}$$
$$\begin{CD} 
[ \cdots \Omega^k_{M'/A} \otimes F_{p+k}V_{\alpha}{\mathcal M} @> \nabla_{/A}^k >>   \Omega^{k+1}_{M'/A} \otimes F_{p+k+1}V_{\alpha}{\mathcal M} \cdots]\\
 @V \:\: \nabla_t^k VV @V  \nabla_t^{+1}k V \:\:V \\
[ \cdots \Omega^1_{A}\otimes \Omega^k_{M'/A}\otimes F_{p+k+1}V_{\alpha+1}{\mathcal M} @> \nabla_{/A}^k >>   \Omega^1_{A} \otimes\Omega^{k+1}_{M'/A}\otimes 
 F_{p+k+2}V_{\alpha+1}{\mathcal M} \cdots] \:.
  \end{CD}$$
$ $\\ 
Again all differentials in the $F$-graded and $V$-filtered complex
\begin{equation}
\begin{CD} gr^F_pV_{\alpha}DR^*_{/A}({\mathcal M})  \\
@V  gr^F(\nabla_t) VV\\
\Omega^1_{A}\otimes gr^F_{p+1}V_{\alpha+1}DR^*_{/A}({\mathcal M})
  \end{CD}
\end{equation}
are ${\mathcal O}_{M'}$-linear. Moreover, the induced vertical ${\mathcal O}_{M'}$-linear  maps
$$gr^F_p(\nabla_t^k): \Omega^k_{M'/A} \otimes 
gr^F_pV_{\alpha}{\mathcal M}\to \Omega^1_{A} \otimes\Omega^k_{M'/A}\otimes  
gr^F_{p+1}V_{\alpha+1}{\mathcal M}; $$
$$\omega \otimes [m] \mapsto  
(-1)^k dt\otimes\omega \otimes (\partial_t [m]) $$
correspond up to a sign $(-1)^k$ under the identification 
$$\Omega^1_{A} \simeq {\mathcal O}_{M'}; \: fdt \simeq f$$
to the maps induced from 
$$\partial_t: (V_{\alpha} {\mathcal M},F_p) \to (V_{\alpha+1} {\mathcal M},F_{p+1})$$
by tensoring with the free (and therefore flat) ${\mathcal O_{M'}}$-modules $\Omega^k_{M'/A}$!\\

So by Lemma 0.1 we get the
\begin{cor} The horizontal inclusion of $F$-filtered double complexes
$$\begin{CD} F_pV_{-1}DR^*_{/A}({\mathcal M}) @>>> F_pDR^*_{/A}({\mathcal M}) \\
@V  \nabla_t VV  @V  \nabla_t VV\\
\Omega^1_{A}\otimes F_{p+1}V_{0}DR^*_{/A}({\mathcal M}) @>>>
\Omega^1_{A}\otimes F_{p+1}DR^*_{/A}({\mathcal M})
  \end{CD}$$
induces a {\em filtered quasi-isomorphism} of the correspondung total complexes, i.e. the horizontal map
$$\begin{CD} gr^F_pV_{-1}DR^*_{/A}({\mathcal M}) @>>> gr^F_pDR^*_{/A}({\mathcal M}) \\
@V  gr^F_p(\nabla_t) VV  @V  gr^F_p(\nabla_t) VV\\
\Omega^1_{A}\otimes gr^F_{p+1}V_{0}DR^*_{/A}({\mathcal M}) @>>> 
\Omega^1_{A}\otimes gr^F_{p+1}DR^*_{/A}({\mathcal M})
  \end{CD}$$
induces a {\em quasi-isomorphism} of the correspondung total complexes for all
$p\in \bb{Z}$. 
\hfill $\Box$
\end{cor}
In fact all vertical maps of the corresponding quotient double complex
$$\begin{CD} gr^F_pDR^*_{/A}({\mathcal M}/V_{-1}{\mathcal M}) \\
@V  gr^F_p(\nabla_t) VV  \\
\Omega^1_{A}\otimes gr^F_{p+1}DR^*_{/A}({\mathcal M}/V_{0}{\mathcal M}) 
  \end{CD}$$
are isomorphisms by  Lemma 0.1 so that the total complex of this quotient double complex is {\em acyclic}! Note that the filtered relative de Rham complex 
$$\left(DR^*_{/A}(-),F\right)$$
is {\em functorial} for filtered $({\mathcal D}_{M'/A},F)$-modules, where the sheaf of relative differential operators ${\mathcal D}_{M'/A}\subset V_0{\mathcal D}_{M'}$ is again filtered
by the order of a differential operator.\\

So the total complex of 
\begin{equation}
\begin{CD} gr^F_pV_{-1}DR^*_{/A}({\mathcal M})  \\
@V  gr^F_p(\nabla_t) VV \\
\Omega^1_{A}\otimes gr^F_{p+1}V_{0}DR^*_{/A}({\mathcal M})
  \end{CD}
  \end{equation}
  represents $gr^F_p DR^*({\mathcal M})$. Therefore it belongs to $D^b_{coh,X'}(M')$, is acyclic for almost all $p$, and can be used for the calculation of 
  $$MHC_y(({\mathcal M},F)) \quad  \text{and} \quad i^! MHC_y(({\mathcal M},F)) \:.$$

In fact for the calculation of $i^! MHC_y(({\mathcal M},F))$ we can even use the total complex of
\begin{equation}
\begin{CD} gr^F_pV_{<-1}DR^*_{/A}({\mathcal M})  \\
@V  gr^F_p(\nabla_t) VV \\
\Omega^1_{A}\otimes gr^F_{p+1}V_{<0}DR^*_{/A}({\mathcal M}) \:,
  \end{CD}
  \end{equation}
because the horizontal complexes of the quotient double complex
$$\begin{CD} gr^F_p gr^V_{-1}DR^*_{/A}({\mathcal M})  \\
@V  gr^F_p(\nabla_t) VV \\
\Omega^1_{A}\otimes gr^F_{p+1} gr^V_{0}DR^*_{/A}({\mathcal M}) \:
  \end{CD}$$
are given by
\begin{equation*}\begin{split}
gr^F_p gr^V_{-1}DR^*_{/A}({\mathcal M}) &\simeq gr^F_p DR^*_{/A}(gr^V_{-1}{\mathcal M})\\
&\simeq i_*(\: gr^F_p DR^*_{/A}(gr^V_{-1}{\mathcal M}|M)\:)
\end{split}
\end{equation*}
and 
\begin{equation*}\begin{split}
gr^F_{p+1} gr^V_{0}DR^*_{/A}({\mathcal M}) &\simeq gr^F_{p+1} DR^*_{/A}(gr^V_{0}{\mathcal M})\\
& \simeq i_*(\: gr^F_{p+1} DR^*_{/A}(gr^V_{0}{\mathcal M}|M)\:)\:,
\end{split}\end{equation*}
if we make the identification $\Omega^1_{A} \simeq {\mathcal O}_{M'}dt$.
Here $i$ is the closed inclusion $M\to M'$, with 
$$F_pDR^*(gr^V_{-1}{\mathcal M}|M) = F_pDR^*(\Psi_{t,-1}({\mathcal M}))$$
and
$$F_{p+1}DR^*(gr^V_{0}{\mathcal M}|M) = F_pDR^*(\Phi_{t,uni}({\mathcal M}))$$
the absolute filtered de Rham complexes of the coherent left ${\mathcal D}_M$-modules
$\Psi_{t,-1}({\mathcal M})$ and $\Phi_{t,uni}({\mathcal M})$
with its induced (shifted) filtrations $F$ (which are good by assumption (s1.)). And from 
$$ gr^F_p DR^*(gr^V_{\alpha}{\mathcal M}|M) \in D^b_{coh,X}(M)$$ 
we get
$$gr^F_p gr^V_{\alpha}DR^*_{/A}({\mathcal M}) \simeq i_*(\: gr^F_p DR^*(gr^V_{\alpha}{\mathcal M}|M)\:)
\in D^b_{coh,X'}(M') \:,$$
together with 
$$i^![gr^F_p gr^V_{\alpha}DR^*_{/A}({\mathcal M})] = i^!i_*[gr^F_p DR^*(gr^V_{\alpha}{\mathcal M|M})]
= 0 \in G_0(X)$$
for $\alpha= -1$ and $0$ by (\ref{i!i*}).\\

So also the total complex of
\begin{equation}
\begin{CD} gr^F_pV_{<-1}DR^*_{/A}({\mathcal M})  \\
@V  gr^F_p(\nabla_t) VV \\
\Omega^1_{A}\otimes gr^F_{p+1}V_{<0}DR^*_{/A}({\mathcal M}) \:,
  \end{CD}
  \end{equation}
  has bounded coherent cohomology supported on $X'$, and can be used for the calculation of
 $i^! MHC_y(({\mathcal M},F))$. The advantage of this double complex is the fact,
 that it is only related to the nearby cycles $\Psi_{t,\alpha}({\mathcal M})$ ($-1\leq \alpha <0$) and not to the unipotent vanishing cycles $\Phi_{t,uni}({\mathcal M})$!\\
  
Next we want to calculate the class
 of its pullback $Li^*(-)$ in $G_0(X)$. But for the horizontal subcomplexes we get by
Lemma 0.2:
\begin{equation}\begin{split}
Li^* gr^F_{p} V_{<-1}DR^*_{/A}( {\mathcal M}) &\simeq 
i^* gr^F_{p} V_{<0} DR^*_{/A}( {\mathcal M})\\ 
&\simeq gr^F_{p} DR(\;V_{<0}/ V_{<-1}{\mathcal M}|M\;) 
\end{split}\end{equation}
and
\begin{equation}\begin{split}
Li^* gr^F_{p+1} V_{<0} DR^*_{/A}( {\mathcal M})&\simeq  
i^* gr^F_{p+1} V_{<0}  DR^*_{/A}({\mathcal M})\\ 
&\simeq gr^F_{p+1} DR(\;V_{<0}/ V_{<-1}{\mathcal M}|M\;) \:.
\end{split}\end{equation}

Putting everything together, we get in $G_0(X)$ by additivity the equality:
\begin{equation}\begin{split}
&i^![gr^F_p DR^*({\mathcal M})] =\\ -[gr^F_{p} DR^*(\;V_{<0}/ V_{<-1}&{\mathcal M}|M\;)]
+ [gr^F_{p+1} DR^*(\;V_{<0}/ V_{<-1}{\mathcal M}|M\;)]  \:.
\end{split}\end{equation}
Note that the $-$-sign for the first class on the right side is coming from the fact,
that $gr^F_{p} DR(\;V_{<0}/ V_{<-1}{\mathcal M}|M\;)$ viewed as a subcomplex of the double complex above agrees only up to a shift by one with the usual convention that ``top-dimesional form'' are in degree zero!\\ 

By using the filtration $V_{\beta}/ V_{<-1}{\mathcal M}$ of
$V_{<0}/ V_{<-1}{\mathcal M}$ by ${\mathcal D}_{M'/A}$-modules
($-1\leq \beta <0$), we get by the assumption (s1.) and additivity
the following equality in the Grothendieck group $ G_0(X)$:
\begin{equation}\begin{split}
&i^![gr^F_p DR^*({\mathcal M})] =\\ \sum_{-1\leq \beta <0}\:(\:
-[gr^F_{p} DR^*(\;&gr^V_{\beta}{\mathcal M}|M\;)]
+ [gr^F_{p+1} DR^*(\;gr^V_{\beta}{\mathcal M}|M\;)]\:) \:.
\end{split}\end{equation}
And this implies Theorem 0.1:
$$i^!MHC_y(\;({\mathcal M},F)\;) =$$
$$ i^!(\;\sum_p\; [gr^F_{-p} DR({\mathcal M})]\cdot (-y)^p\;) \: =$$
$$\sum_{p,-1\leq \beta <0}\:(\: -[gr^F_{-p} DR(\Psi_{t,\beta}{\mathcal M})]\cdot (-y)^p
+ [gr^F_{-(p-1)} DR(\Psi_{t,\beta}{\mathcal M})]\cdot (-y)^p \:) $$
$$= -(1+y)\cdot \sum_{-1\leq \beta <0}\:MHC_y(\;\Psi_{t,\beta}({\mathcal M},F)\;)$$
$$= -(1+y)\cdot MHC_y(\;\Psi_{t}({\mathcal M},F)\;) \:.$$ \hfill $\Box$

\section{Specialization for Mixed Hodge modules}
We can apply our first main theorem \ref{mainthm} in the context of algebraic mixed Hodge modules
in the following more general context (for a base field $k\subset \bb{C}$). Let $f: X'\to \bb{A}^1_k$ be an algebraic function on the algebraic variety $X'$ over $k$ such that the inclusion of the zero fiber
$i: X:=\{f=0\}\to X'$ is everywhere of codimension one (i.e. $f$ is not vanishing on
any irreducible component of $X'$). Then one gets as before a Gysin map
$i^!: G_0(X')\to G_0(X)$ together with an exact nearby cycle transformation
$\Psi_f: MHM(X'/k) \to  MHM(X/k)$ on the abelian category of algebraic mixed Hodge modules.
This induces therefore also a transformation of the corresponding Grothendieck groups
$$\Psi^H_f: K_0(MHM(X'/k)) \to  K_0(MHM(X/k)) \:.$$
In the introduction we already explained that by \cite{BSY} we have in this context a motivic Hodge-Chern class transformation 
$$MHC_y: K_0(MHM(-)) \to G_0(-)[y,y^{-1}] \:,$$
which {\em commutes with proper push down}. Let us recall our main result in this context.

\begin{thm} 
This motivic Hodge-Chern class transformation commutes with specialization in the following sense:
\begin{equation}
-(1+y)\cdot MHC_y(\; \Psi^H_f (-)\;) = i^!MHC_y(-)   
\end{equation} 
as transformations $K_0(MHM(X'/k))\to G_0(X)[y,y^{-1}]$.
\end{thm} 

Its proof can be reduced to theorem \ref{mainthm}, since all transformations $MHC_y$,
$\Psi^H_f$ and $i^!$ commute with proper pushdown. 
By {\em resolution of singularities} and ``additivity'', $G_0(X')$ is generated by
classes $\pi_*[{\mathcal M}]$ of mixed Hodge modules ${\mathcal M}$ on an algebraic manifold
$M$, with $\pi: M\to X'$ a proper morphism. It is enough to prove the stated result for such
a generating class $\pi_*[{\mathcal M}]$. If a connected component $S$ of $M$ maps into
$X$, then this doesn't contribute to the specialization:
$$  \Psi^H_f (\pi_*[{\mathcal M}|S]) = 0 \quad \text{and} \quad i^!MHC_y(\pi_*[{\mathcal M}|S]) = 0 \:.$$
So we can assume $f':=f\circ \pi$ is not vanishing on any connected component of $M$.
Since all involved transformations $MHC_y, \Psi^H_f$
and $i^!$ commute with $\pi_*$, we can assume $X'=M$ is smooth (with $\pi_*=id_*$). 
Using ``additivity'' and the weight filtration, one can even reduce to the case of pure Hodge modules. 
By using the (proper) graph embedding of $f'$, we can reduce to the case $M':=M\times \bb{A}^1_k$ and $t$ the projection onto $\bb{A}^1_k$ as studied before (with $X':=graph_f(M)\subset M'$).
Moreover the nearby cycles $\Psi^H_t$ of the filtered ${\mathcal D}$-module
$({\mathcal M},F)$  underlying a pure Hodge module on $M'$ is exactly given as in definition \ref{psi}.
And  $({\mathcal M},F)$ is by Saito \cite{Sa1} strictly specializable (as explained in remark
\ref{pure}). So in this case the claim is a special case of theorem \ref{mainthm}. \hfill $\Box$\\

In the complex analytic context one doesn't have the full calculus of Grothendieck functors
on the derived category of mixed Hodge modules. But nevertheless in the embedded context
of an analytic subset $X\subset M$ in a complex manifold $M$, with $f: M\to \bb{C}$
a holomorphic function such that $X_0:=X\cap \{f=0\}$ is {\em compact}, one has an {\em exact}
nearby cycle functor 
$$\Psi^H_f: MHM(X)\simeq MHM_X(M) \to MHM_{X_0}(M)\simeq  MHM(X_0)\:,$$
with $MHM_X(M)$ the abelian category of mixed Hodge modules on $M$ with support in $X$, and similarly for $MHM_{X_0}(M)$ (compare \cite{Sa2}).
Also the graph embedding 
$$g:=(id_M,f): M\to M':=M\times \bb{C} \quad \text{with} \quad
g: X\simeq X':=g(X)$$
induces an exact functor $g_*: MHM_X(M) \to MHM_{X'}(M')$ commuting with 
$\Psi^H_f$ and $\Psi^H_t$. 
Similarly one gets a motivic Hodge Chern class transformation
$$MHC_y: K_0(MHM(X))=K_0(MHM_X(M))\to G^X_0(M)[y,y^{-1}] \:.$$
Here one uses the fact that morphisms of analytic mixed Hodge modules
are {\em strict} with respect to the (Hodge) filtration $F$ of the underlying
filtered ${\mathcal D}$-modules on $M$.  And also $MHC_y$ commutes with
$g_*: G^X_0(M)[y,y^{-1}]\to G^{X'}_0(M')[y,y^{-1}]$.\\
 
So in this embedded complex analytic context we can argue as
in the end of the argument before, and get the

\begin{thm} Let $f: M\to \bb{C}$ be a holomorphic function on the complex manifold $M$,
with $X\subset M$ a closed analytic subset. Assume $\{f=0\}$ is a hypersurface of codimension one, with $X_0:=X\cap \{f=0\}$ compact. Then the motivic Hodge-Chern class 
transformation $MHC_y$ commutes with specialization in the following sense:
\begin{equation}
-(1+y)\cdot MHC_y(\; \Psi^H_f (-)\;) = i^!MHC_y(-)   
\end{equation} 
as transformations 
$$K_0(MHM(X))=K_0(MHM_X(M))\to G^{X_0}_0(M)[y,y^{-1}]\:.$$\hfill $\Box$
\end{thm}

\end{document}